\documentclass{amsart}
\usepackage{amsfonts}

\newtheorem{theorem}{Theorem}
\theoremstyle{plain}

\newtheorem{corollary}{Corollary}

\newtheorem{lemma}{Lemma}

\newtheorem{proposition}{Proposition}

\numberwithin{equation}{section}

\begin{document}
\title{On restriction properties of equivariant K-theory rings.}
\author{Abdelouahab AROUCHE}
\address{USTHB, Dept. Math. BP 32, El Alia 16111 Alger(ie)}
\date{December 21st, 2004}
\maketitle

\begin{abstract}
An important ingredient in the completion theorem of equivariant K-theory
given by S. Jackowski is that the representation ring R($\Gamma $) of a
compact Lie group satisfies two restriction properties called ($N$) and ($R_{%
\mathcal{F}}$).

We give in this note sufficient conditions on a (compact) $\Gamma $-space $Z$
such that these properties hold with $K_{\Gamma }^{*}(Z)$ instead of $%
R(\Gamma )$. As an example, we consider the space $Z(\Gamma ;G)$ of the so
called ``elementary cocycles with coefficients in $G$'' invented by H.
Ibisch in his construction of a universal $(\Gamma ;G)$-bundle.
\end{abstract}

\section{Introduction}

In what follows, let $\Gamma $ be a compact Lie group and $Z$ a compact $%
\Gamma $-space. Consider, for every subgroup $\Lambda \leq \Gamma $, the
ring $M_{\Lambda }=K_{\Lambda }^{*}(Z)$. We have a restriction morphism $%
res:M_{\Gamma }\rightarrow M_{\Lambda }$ whose kernel is denoted by $%
I_{\Lambda }^{\Gamma }$ or $I_{\Lambda }$ for short. If $\mathcal{F}$ is a
family of (closed) subgroups of $\Gamma $ (stable under subconjugation),
then the $\mathcal{F}$-adic topology on $M_{\Gamma }$ is defined by the set
of ideals
\begin{equation*}
I(\mathcal{F})=\left\{ I_{\Lambda _{1}}.I_{\Lambda _{2}}...I_{\Lambda
_{n}},\ \Lambda _{i}\in \mathcal{F},\ i=1,...,n\right\} .
\end{equation*}

If $\Lambda \leq \Gamma $, then $\mathcal{F}\cap \Lambda =\left\{ \Omega \in
\mathcal{F}:\Omega \leq \Lambda \right\} $ is a family of subgroups of $%
\Lambda $.

\bigskip

Let $X$ be a compact $\Gamma $-space with a $\Gamma $-map $\sigma
:X\rightarrow Z$. Define the two pro-rings
\begin{equation*}
K_{\Gamma }^{*}\left[ \mathcal{F}\right] (X)=\left\{ K_{\Gamma }^{*}\left(
X\times K\right) :\ K\ compact\ \subseteq E\mathcal{F}\right\} ,
\end{equation*}

and
\begin{equation*}
\left( K_{\Gamma }^{*}/I\left( \mathcal{F}\right) \right) \left( X\right)
=\left\{ K_{\Gamma }^{*}\left( X\right) /I.K_{\Gamma }^{*}\left( X\right) :\
I\in I\left( \mathcal{F}\right) \right\}
\end{equation*}

Here $E\mathcal{F}$ denotes the classifying space of T. tom Dieck [7]. The
completion theorem asserts that the pro-homomorphism
\begin{equation*}
p_{X}(\mathcal{F}):K_{\Gamma }^{*}\left[ \mathcal{F}\right] (X)\rightarrow
\left( K_{\Gamma }^{*}/I\left( \mathcal{F}\right) \right) \left( X\right)
\end{equation*}

induced by the projection $X\times E\mathcal{F}\rightarrow X$, is actually
an isomorphism, provided $M_{\Lambda }$ satisfies the conditions $(N)$ and $%
(R_{\mathcal{F}})$, and $K_{\Lambda }(X)$ is finitely generated over $%
M_{\Lambda }$ (via $\sigma $), for all $\Lambda \leq \Gamma $.

\bigskip

In order to prove his completion theorem [4], S. Jackowski required the
following two conditions, and showed they are satisfied when $Z=*$ is a
point :

$(N)$ : $M_{\Gamma }$ is noetherian and $M_{\Lambda }$ is a finitely
generated $M_{\Gamma }$-module, $\forall \Lambda \leq \Gamma $.

$(R_{\mathcal{F}})$ : For every family $\mathcal{F}$ of subgroups of $\Gamma
$, and every subgroup $\Lambda \leq \Gamma $, the $\mathcal{F}$-adic
topology on $M_{\Lambda }$ defined by the restriction $res:M_{\Gamma
}\rightarrow M_{\Lambda }$ coincides with the $\left( \mathcal{F}\cap
\Lambda \right) $-adic topology.

\section{Some nilpotent elements in $K_{\Gamma }^{*}(Z).$}

Let $z\in Z$, and $\Gamma _{z}$ be its isotropy subgroup. Then an element $%
\zeta \in M_{\Gamma }$ belongs to $\ker \left\{ K_{\Gamma }^{*}(Z)\overset{%
\varphi _{z}}{\rightarrow }R(\Gamma _{z})\right\} $ if and only if its
restriction to the orbit $\Gamma .z$ is zero. It is well known that the
elements of $\ker \left\{ K_{\Gamma }^{*}(Z)\rightarrow \underset{z\in Z}{%
\prod }R(\Gamma _{z})\right\} $ are nilpotent [6](5.1). We have the
following :

\begin{lemma}
Let $z\in Z$ and $\zeta \in \ker \left\{ K_{\Gamma }^{*}(Z)\rightarrow
R(\Gamma _{z})\right\} $. Then there exists a relatively compact open $%
\Gamma $-neighbourhood $V(\zeta )$ of $z$ such that the restriction of $%
\zeta $ to $\overline{V(\zeta )}$ is nilpotent.
\end{lemma}

\begin{proof}
If $\zeta =E-\mathbf{T}\in \ker \left\{ K_{\Gamma }^{*}(Z)\rightarrow
R(\Gamma _{z})\right\} $, then $E_{z}-T=0\in R(\Gamma _{z})$, where $T$ is
seen as a $\Gamma _{z}$-module. Therefore, there exists a $\Gamma _{z}$%
-module $R$ such that $E_{z}\oplus R\cong T\oplus R$. Consider the
isomorphism between $\Gamma _{z}$- vector bundles :
\begin{equation*}
E\oplus \mathbf{R\mid }_{\left\{ z\right\} }\mathbf{\cong T\oplus R\mid }%
_{\left\{ z\right\} }.
\end{equation*}

According to [6](1.2), there is a $\Gamma _{z}$-neighbourhood $U$ of $z$ and
a $\Gamma _{z}$- isomorphism
\begin{equation*}
E\oplus \mathbf{R\mid }_{U}\mathbf{\cong T\oplus R\mid }_{U}.
\end{equation*}

Let $S$ be a slice at $z$ ([2] II.4.1). Then $\Gamma .\left( S\cap U\right) $
is an open $\Gamma $-neighbourhood of $z$ such that $\zeta \in \ker \left\{
K_{\Gamma }^{*}(Z)\rightarrow R(\Gamma _{y})\right\} $, $\forall y\in \Gamma
.\left( S\cap U\right) $. Now let $V(\zeta )$ be any open $\Gamma $%
-neighbourhood of $z$ satisfying
\begin{equation*}
z\in V(\zeta )\subseteq \overline{V(\zeta )}\subseteq \Gamma .\left( S\cap
U\right) .
\end{equation*}
\end{proof}

\begin{proposition}
Assume $K_{\Gamma }^{*}(Z)$ is noetherian. Then there exists a finite subset
$Y\subseteq Z$ such that the ideal $\ker \left\{ K_{\Gamma
}^{*}(Z)\rightarrow \underset{y\in Y}{\prod }R(\Gamma _{y})\right\} $
consists of nilpotent elements.
\end{proposition}

\begin{proof}
For each $z\in Z$, let $\ker \left\{ K_{\Gamma }^{*}(Z)\rightarrow R(\Gamma
_{z})\right\} $ be generated by $\zeta _{1},...,\zeta _{n_{z}}$, and for $%
j=1,...,n_{z}$, let $V(\zeta _{j})$ be as in lemma (2.1). Put
\begin{equation*}
V_{z}=\underset{j=1,...,n_{z}}{\cap }V(\zeta _{j}).
\end{equation*}

Since $Z$ is compact, it can be covered by a finite number of such open
subsets :
\begin{equation*}
Z=\underset{i=1,...,n}{\cup }V_{z_{i}}.
\end{equation*}

It is easy to show that the ideal $\ker \left\{ K_{\Gamma
}^{*}(Z)\rightarrow \underset{i=1,...,n}{\prod }R(\Gamma _{z_{i}})\right\} $
consists of nilpotent elements.
\end{proof}

\section{Supports of primes}

In [4], G. Segal defined the support \textit{supp}$y$ of a prime $y\in
SpecR\left( \Gamma \right) $ to be a minimal subgroup from which $y$ comes.
We have shown in [1] that such a notion does exist for primes in $M_{\Gamma
} $. Indeed, if $p:Z\rightarrow *$ is the projection of a $\Gamma $-space $Z$
on a fixed point and $p^{*}:R\left( \Gamma \right) \rightarrow M_{\Gamma }$
is the induced homomorphism, we have the following :

\begin{proposition}
Let $\Lambda \leq \Gamma $ be a (closed) subgroup. Assume $M_{\Gamma }$ is
noetherian and let $x\in SpecM_{\Gamma }$. Then the following properties are
equivalent :

\begin{enumerate}
\item $x$ come from $M_{\Lambda }$, i.e. $\exists y\in SpecM_{\Lambda
}:x=res^{-1}\left( y\right) $.

\item $x$ contains $I_{\Lambda }$.

\item \textit{supp}$\left( p^{*-1}\left( x\right) \right) $ is subconjugate
to $\Lambda $.
\end{enumerate}
\end{proposition}

In order to prove proposition (3.1), we need the following lemma :

\begin{lemma}
Let $x\in SpecM_{\Gamma }$ and $z\in Z$ be such that there exists a prime $%
t\in SpecR\left( \Gamma _{z}\right) $, with support $\Sigma $, satisfying $%
x=\varphi _{z}^{-1}\left( t\right) ,\ \varphi _{z}:M_{\Gamma }\rightarrow
R\left( \Gamma _{z}\right) .$ Let $\gamma \in \Gamma $ and $z^{\prime }$=$%
\gamma .z$. Then there exists a prime $t^{\prime }\in SpecR\left( \Gamma
_{z^{\prime }}\right) $, with support $\Sigma ^{\prime }=\gamma \Sigma
\gamma ^{-1}$, satisfying $x=\varphi _{z^{\prime }}^{-1}\left( t^{\prime
}\right) ,\ \varphi _{z^{\prime }}:M_{\Gamma }\rightarrow R\left( \Gamma
_{z^{\prime }}\right) .$
\end{lemma}

\begin{proof}
Let $\Gamma .z\overset{i}{\hookrightarrow }Z$ be the inclusion of the orbit
of $z$. We have the following commutative diagram :
\begin{equation*}
\begin{array}{ccc}
& R\left( \Gamma _{z}\right) &  \\
\overset{\varphi _{z}}{\nearrow } &  & \overset{\chi _{z}}{\nwarrow } \\
K_{\Gamma }^{*}\left( Z\right) & \overset{i^{*}}{\rightarrow } & K_{\Gamma
}^{*}\left( \Gamma .z\right) \\
\underset{\varphi _{\gamma .z}}{\searrow } &  & \underset{\chi _{\gamma .z}}{%
\swarrow } \\
& R\left( \Gamma _{\gamma .z}\right) &
\end{array}%
\end{equation*}

It is clear that both $\chi _{z}$ and $\chi _{\gamma .z}$ are isomorphisms.
Put $h=\chi _{z}^{-1}\circ \chi _{\gamma .z}$. Then $h$ admits a restriction
$\widetilde{h}:R\left( \Sigma \right) \rightarrow R\left( \gamma \Sigma
\gamma ^{-1}\right) $, which makes the following diagram commute and ends
the proof :

\begin{equation*}
\begin{array}{ccccccc}
&  &  &  & R\left( \Gamma _{z}\right) & \overset{res}{\rightarrow } &
R\left( \Sigma \right) \\
&  &  & \overset{\chi _{z}}{\nearrow } &  &  &  \\
K_{\Gamma }^{*}\left( Z\right) & \overset{i^{*}}{\rightarrow } & K_{\Gamma
}^{*}\left( \Gamma .z\right) &  & \cong \downarrow h &  & \cong \downarrow
\widetilde{h} \\
&  &  & \underset{\chi _{\gamma .z}}{\searrow } &  &  &  \\
&  &  &  & R\left( \Gamma _{\gamma .z}\right) & \overset{res}{\rightarrow }
& R\left( \gamma \Sigma \gamma ^{-1}\right)%
\end{array}%
\end{equation*}
\end{proof}

\bigskip

\begin{proof}
Now we proceed to prove proposition (3.1) (see [1]). $1)\Longrightarrow 2)$
and $2)\Longrightarrow 3)$ are obvious. Let us prove $3)\Longrightarrow 1)$.
To this end, let $z_{1},...,z_{n}$ be elements of $Z$ such that $\ker
\left\{ \varphi :K_{\Gamma }^{*}(Z)\rightarrow \underset{i=1,...,n}{\prod }%
R(\Gamma _{z_{i}})\right\} $ consists of nilpotent elements. Since every $%
R(\Gamma _{z_{i}})$ is finitely generated over $M_{\Gamma }$ (it i even
finitely generated over $R\left( \Gamma \right) $ [5]), the induced map on
spectra $\varphi _{*}:\underset{i=1,...,n}{\coprod }SpecR(\Gamma
_{z_{i}})\rightarrow SpecM_{\Gamma }$ is onto. So there is an element $z\in
Z $ and a prime $t\in SpecR(\Gamma _{z})$ such that $x=\varphi _{z}^{-1}(t)$%
. Following [5] (3.7), we have s\textit{upp}$(p^{*-1}(x))$=\textit{supp}$%
\left( t\right) $. Let $\Sigma =\mathit{supp}(p^{*-1}(x))$. By assumption,
the exists $\gamma \in \Gamma $ such that $\gamma \Sigma \gamma ^{-1}\leq
\Lambda $. Put $z^{\prime }=\gamma .z$. According to lemma (3.2), there
exists a prime $t^{\prime }\in SpecR(\Gamma _{z^{\prime }})$, with support $%
\Sigma ^{\prime }=\gamma \Sigma \gamma ^{-1}$, and satisfying
\begin{equation*}
x=\varphi _{z^{\prime }}^{-1}\left( t^{\prime }\right) ,\ \varphi
_{z^{\prime }}:M_{\Gamma }\rightarrow R\left( \Gamma _{\gamma z}\right) .
\end{equation*}

Therefore, $x\in SpecM_{\Gamma }$ comes from $\Lambda $, as shown by the
commutative diagram :
\begin{equation*}
\begin{array}{ccc}
M_{\Gamma } & \rightarrow & M_{\Lambda } \\
\downarrow &  & \downarrow \\
R(\Gamma _{\gamma .z}) & \rightarrow & R\left( \Lambda _{\gamma .z}\right)
\\
\searrow &  & \ \swarrow \\
& R\left( \Sigma ^{\prime }\right) &
\end{array}%
\end{equation*}
\end{proof}

\bigskip

Hence, \textit{supp}$\left( p^{*-1}\left( x\right) \right) $ is, up to
conjugation, the minimal subgroup of $\Gamma $ from which $x$ comes. We
denote it by \textit{supp}$\left( x\right) $. Consequently, if $\left(
x,y\right) \in SpecM_{\Gamma }\times SpecM_{\Lambda }$ is such that $%
x=res^{-1}\left( y\right) $, then \textit{supp}$\left( x\right) =$\textit{%
supp}$(y)$.

\section{The restriction properties}

Let us first consider the condition $(R_{\mathcal{F}})$. We have the
following :

\begin{theorem}
Assume $M_{\Lambda }$ is noetherian $\forall \Lambda \leq \Gamma $. Then,
for every family $\mathcal{F}$ of subgroups of $\Gamma $ , the condition $%
(R_{\mathcal{F}})$ is satisfied.
\end{theorem}

\begin{proof}
For each $\Omega \in \mathcal{F}$ we have $I_{\Omega }^{\Lambda }\supseteq
I_{\Omega }^{\Gamma }.M_{\Lambda }$. So it remains only to show that for any
ideal $I.M_{\Lambda }$, where $I\in I\left( \mathcal{F}\right) $, there is
an ideal $J\in I\left( \mathcal{F}\cap \Lambda \right) $ such that $%
J\subseteq r\left( I.M_{\Lambda }\right) $, since $M_{\Lambda }$ is
noetherian. It is enough to do that for $I=I_{\Omega }^{\Gamma }$. So let $%
p_{1},...,p_{n}$ be minimal prime ideals in $M_{\Lambda }$ containing $%
I.M_{\Lambda }$, and let $\Sigma _{1},...,\Sigma _{n}$ be their supports.
Put $J=I_{\Sigma _{1}}...I_{\Sigma _{n}}$. Since $p_{i}$ comes from $\Sigma
_{i}$, we have $J\subseteq I_{\Sigma _{1}}\cap ...\cap I_{\Sigma
_{n}}\subseteq p_{1}\cap ...\cap p_{n}=r\left( I\right) $. Now, if we put $%
q_{i}=res^{-1}\left( p_{i}\right) $, then $\Sigma _{i}=$\textit{supp}$%
(q_{i}) $. But $q_{i}$ comes from $\Omega $ because it contains $I_{\Omega
}^{\Gamma }$. Hence, $\Sigma _{i}$ is subconjugate to $\Omega $. Since $%
\Omega \in \mathcal{F}$, so doe $\Sigma _{i}$, for $i=1,...,n$, that is, $%
J\in I\left( \mathcal{F}\cap \Lambda \right) $.

\begin{corollary}
If the condition $(N)$ is satisfied, then so is the condition $(R_{\mathcal{F%
}})$.
\end{corollary}
\end{proof}

\bigskip

It is easy to see that the condition $\left( N\right) $ is fulfilled if $Z$
is a (compact) differentiable $\Gamma $-manifold, and more generally, for
any finite $\Gamma $-CW-complex, for $M_{\Lambda }$ is a finitely generated
module over $R\left( \Lambda \right) $, $\forall \Lambda \leq \Gamma $. This
can be shown by induction and use of Mayer-Vietoris sequence.

\begin{corollary}
Assume $M_{\Lambda }$ is finitely generated over $R\left( \Gamma \right) $, $%
\forall $ $\Lambda \leq \Gamma $. Then the conditions $\left( N\right) $ and
$\left( R_{\mathcal{F}}\right) $ are satisfied for every family $\mathcal{F}$
of subgroups of $\Gamma $.
\end{corollary}

\section{An example}

For a topological group $G$ and a family $\mathcal{F}$ of subgroups of a
compact Lie group $\Gamma $, H.Ibisch has constructed in [2] a universal $%
\left( \Gamma ;G\right) $-bundle $E\left( \mathcal{F};G\right) $ whose
``components'' are the spaces
\begin{equation*}
Z\left( \Lambda ;G\right) =\left\{ f:\Lambda \times \Lambda \rightarrow
G:f\left( \lambda _{1}.\lambda _{2},\lambda _{3}\right) =f\left( \lambda
_{1},\lambda _{2}.\lambda _{3}\right) ,\ \forall \lambda _{1},\lambda
_{2},\lambda _{3}\in \Lambda \right\} ,\Lambda \leq \Gamma .
\end{equation*}

The crucial fact is that for a $\Gamma $-space $X$, trivial $\left( \Gamma
;G\right) $-bundles over $X$ correpond to $\Gamma $-maps $\sigma
:X\rightarrow Z\left( \Gamma ;G\right) $. We then take $Z=Z\left( \Gamma
;G\right) $. Actually, the $\Gamma $-space $Z\left( \Gamma ;G\right) $ turns
to be $\Gamma $-homeomorphic to the space of maps $\Gamma \rightarrow G$
sending $1_{\Gamma }$ to $1_{G}$. Moreover, $Z\left( \Gamma ;G\right)
^{\Gamma }$ is the space of (continuous) homomorphisms $\Gamma \rightarrow G$%
. When $\Gamma $ is finite, the space $Z\left( \Gamma ;G\right) $ has a
further description. For instance, if $\Gamma =\mathbf{Z}_{n}$ is the cyclic
group of order $n$, then $Z\left( \Gamma ;G\right) $ is the product $%
G^{n-1}\times \left\{ 1_{G}\right\} $ with the action $\gamma ^{i}.\left(
g_{1},...,g_{n}=1_{G}\right) =\left(
g_{1+i}g_{i}^{-1},...,g_{n+i}g_{i}^{-1}\right) $, where $\gamma $ is a
generator of $\mathbf{Z}_{n}$ and the indexation is mod $\left[ n\right] $,
hence $Z(\mathbf{Z}_{2};G)$ is no but $G$ with involution $g\mapsto g^{-1}$.
If $G$ is moreover a compact Lie group, then $Z\left( \Gamma ;G\right) $ is
a compact $\Gamma $-manifold. Accordingly, the conditions $\left( N\right) $
and $R\left( \mathcal{F}\right) $ are satisfied.

\end{document}